\documentclass[11pt]{article}

\RequirePackage{amsmath,amsthm,amssymb}

\usepackage[longnamesfirst,sort]{natbib}

\usepackage{geometry}
\usepackage[dvips]{changebar}

\newtheorem{theorem}{Theorem}
\newtheorem{corollary}{Corollary}

\newtheorem{proposition}{Proposition}
\theoremstyle{remark}
\newtheorem{remark}{Remark}
\newtheorem{example}{Example}
\newtheorem{definition}{Definition}

\begin{document}

\title{On the Markov Chain Central Limit Theorem}

\author{Galin L. Jones \\ School of Statistics \\ University of
  Minnesota \\ Minneapolis, MN, USA \\ {\tt galin@stat.umn.edu}}

\maketitle
\begin{abstract}
  The goal of this expository paper is to describe conditions which
  guarantee a central limit theorem for functionals of general state
  space Markov chains.  This is done with a view towards Markov chain
  Monte Carlo settings and hence the focus is on the connections
  between drift and mixing conditions and their implications.  In
  particular, we consider three commonly cited central limit theorems
  and discuss their relationship to classical results for mixing
  processes.  Several motivating examples are given which range from
  toy one-dimensional settings to complicated settings encountered in
  Markov chain Monte Carlo.
\end{abstract}

\section{Introduction}
\label{sec:intro}
Let $X = \{X_i: i=0,1,2,\dots\}$ be a Harris ergodic Markov chain on a
general space $\mathsf{X}$ with invariant probability distribution
$\pi$ having support $\mathsf{X}$.  Let $f$ be a Borel function and
define $\bar{f}_{n} := n^{-1} \sum_{i=1}^n f(X_i)$ and $\text{E}_{\pi}
f : = \int_{\mathsf{X}} f(x) \pi(dx)$.  When $\text{E}_{\pi} |f| <
\infty$ the ergodic theorem guarantees that $\bar{f}_{n} \rightarrow
\text{E}_{\pi} f$ with probability 1 as $n \rightarrow \infty$. The
main goal here is to describe conditions on $X$ and $f$ under which a
central limit theorem (CLT) holds for $\bar{f}_{n}$; that is,
\begin{equation}
\label{eq:clt}
\sqrt{n} (\bar{f}_{n} - \text{E}_{\pi} f) \stackrel{d}{\rightarrow}
\text{N} (0, \sigma^{2}_{f}) 
\end{equation}
as $n \rightarrow \infty$ where $\sigma^{2}_{f} := \text{var}_{\pi} \{
f(X_{0})\} + 2 \sum_{i=1}^{\infty} \text{cov}_{\pi} \{ f(X_{0}),
f(X_{i})\} < \infty$.  Although all of the results presented in this
paper hold in general, the primary motivation is found in Markov chain
Monte Carlo (MCMC) settings where the existence of a CLT is an
extremely important practical problem.  Often $\pi$ is high
dimensional or known only up to a normalizing constant but the value
of $\text{E}_{\pi} f$ is required.  If $X$ can be simulated then
$\bar{f}_{n}$ is a natural estimate of $\text{E}_{\pi} f$.  The
existence of a CLT then allows one to estimate $ \sigma^{2}_{f}$ in
order to decide if $\bar{f}_{n}$ is a good estimate of $\text{E}_{\pi}
f$. (Estimation of $\sigma^{2}_{f}$ is challenging and requires
specialized techniques that will not be considered further here; see
\citet{jone:hara:caff:2004} and \citet{geye:1992} for an
introduction.)  Thus the existence of a CLT is crucial to sensible
implementation of MCMC; see \citet{jone:hobe:2001} for more on this
point of view.  The following simple example illustrates one of the
situations common in MCMC settings.

\begin{example}  \label{ex:hardcore}
  Consider a simple hard-shell (also known as hard-core) model.
  Suppose ${\cal X}= \{ 1, \ldots, n_{1} \} \times \{ 1, \ldots, n_{2}
  \} \subseteq \mathbb{Z}^{2}$.  A \textit{proper} configuration on
  ${\cal X}$ consists of coloring each point either black or white in
  such a way that no two adjacent points are white.  Let $\mathsf{X}$
  denote the set of all proper configurations on ${\cal X}$,
  $N_{\mathsf{X}} (n_{1}, n_{2})$ be the total number of proper
  configurations and $\pi$ be the uniform distribution on $\mathsf{X}$
  so that each proper configuration is equally likely.  Suppose our
  goal is to calculate the typical number of white points in a proper
  configuration; that is, if $W(x)$ is the number of white points in
  $x \in \mathsf{X}$ then we want the value of
\[
\text{E}_{\pi} W = \sum_{x \in \mathsf{X}} \frac{w(x)}{N_{\mathsf{X}}
  (n_{1}, n_{2})} \; .
\]
If $n_{1}$ and $n_{2}$ are even moderately large then we will have to
resort to an approximation to $\text{E}_{\pi} W$.  Consider the
following Markov chain on $\mathsf{X}$.  Fix $p \in (0,1)$ and set
$X_{0} = x_{0}$ where $x_{0} \in \mathsf{X}$ is an arbitrary proper
configuration.  Randomly choose a point $(x,y) \in {\cal X}$ and
independently draw $U \sim \text{Uniform}(0,1)$.  If $u \le p$ and all
of the adjacent points are black then color $(x,y)$ white leaving all
other points alone.  Otherwise, color $(x,y)$ black and leave all
other points alone.  Call the resulting configuration $X_{1}$.
Continuing in this fashion yields a Harris ergodic Markov chain
$\{X_{0}, X_{1}, X_{2}, \ldots\}$ having $\pi$ as its invariant
distribution.  It is now a simple matter to estimate $\text{E}_{\pi}
W$ with $\bar{w}_{n}$.  Also, since $\mathsf{X}$ is finite (albeit
potentially large) it is well known that $X$ will converge
exponentially fast to $\pi$ which implies that a CLT holds for
$\bar{w}_{n}$.
\end{example}

Following the publication of the influential book by
\citet{meyn:twee:1993} the use of drift and minorization conditions
has become a popular method for establishing the existence of a CLT.
Indeed without this constructive methodology it is difficult to
envision how one would deal with complicated situations encountered in
MCMC.  In turn, this has led much of the recent work on general state
space Markov chains to focus on the implications of drift and
minorization.  Another outcome of this approach is that classical
results in mixing processes have been somewhat neglected.  For
example, \citet{numm:2002} and \citet{robe:rose:2004} recently
provided nice reviews of Markov chain theory and its connection to
MCMC.  In particular, both articles contain a review of CLTs for
Markov chains but neither contains any substantive discussion of the
results from mixing processes.  On the other hand, work on mixing
processes rarely discusses their applicability to the important Markov
chain setting outside of the occasional discrete state space example.
For example, \citet{brad:1999} provided a recommended review of CLTs
for mixing processes but made no mention of their connections with
Markov chains.  Also, \citet{robe:1995} gave a brief discussion of the
implication of mixing conditions for Markov chain CLTs but failed to
connect them to the use of drift conditions.  Thus one of the main
goals of this article is to consider the connections between drift and
minorization and mixing conditions and their implications for the CLT
for general state space Markov chains.

\section{Markov Chains and Examples}
\label{sec:markov}
Let $P(x,dy)$ be a Markov transition kernel on a general space
$(\mathsf{X},{\cal B}(\mathsf{X}))$ and write the associated discrete
time Markov chain as $X = \{X_i: i=0,1,2,\dots\}$.  For $n \in
\mathbb{N} := \{1,2,3,\ldots \}$, let $P^n(x,dy)$ denote the $n$-step
Markov transition kernel corresponding to $P$.  Then for $i \in
\mathbb{N}$, $x \in \mathsf{X}$ and a measurable set $A$, $P^{n}(x,A)
= \Pr \left(X_{n+i} \in A|X_i = x\right)$.  Let $f : \mathbb{R}
\rightarrow \mathbb{R}$ be a Borel function and define $P f(x) := \int
f(y) P(x,dy)$ and $\Delta f(x) := Pf(x) - f(x)$.  Always, $X$ will
be assumed to be Harris ergodic, that is, aperiodic,
$\psi$-irreducible and positive Harris recurrent; for definitions see
\citet{meyn:twee:1993} or \citet{numm:1984}. These assumptions are
more than enough to guarantee a strong form of convergence: for every
initial probability measure $\lambda(\cdot)$ on ${\cal B}(\mathsf{X})$
\begin{equation}
\label{eq:bas_conv}
\| P^{n} (\lambda, \cdot) - \pi(\cdot) \| \rightarrow 0 \hspace*{4mm}
\text{ as } \hspace*{3mm} n \rightarrow \infty
\end{equation}
where $ P^{n} (\lambda, A) := \int_{\mathsf{X}} P^{n}(x,A)
\lambda(dx)$ and $\| \cdot \|$ is the total variation norm.
Throughout we will be concerned with the rate of this convergence.
Let $M(x)$ be a nonnegative function and $\gamma(n)$ be a nonnegative
decreasing function on $\mathbb{Z}_{+}$ such that
\begin{equation}
\label{eq:tvbd}
\| P^{n} (x, \cdot) - \pi(\cdot) \| \le M(x) \gamma(n) \; .
\end{equation}
When $X$ is \textit{geometrically ergodic} \eqref{eq:tvbd} holds with
$\gamma(n) = t^{n}$ for some $t < 1$.  \textit{Uniform ergodicity}
means $M$ is bounded and $\gamma(n) = t^{n}$ for some $t < 1$.
\textit{Polynomial ergodicity of order m} where $m \ge 0$ corresponds
to $\gamma(n)=n^{-m}$.

Establishing \eqref{eq:tvbd} directly may be difficult when
$\mathsf{X}$ is a general space.  However, some constructive methods
are given in the following brief discussion; the interested reader
should consult \citet{jarn:robe:2002} and \citet{meyn:twee:1993} for a
more complete introduction to these methods.

A \textit{minorization condition} holds on a set $C$ if there exists a
probability measure $Q$ on ${\cal B}(\mathsf{X})$, a positive integer
$n_0$ and an $\epsilon >0$ such that
\begin{equation}
\label{eq:min_con}
P^{n_0} (x,A) \ge \epsilon \,Q(A) \hspace{3mm} \forall \, x \in C\; ,
\; A \in {\cal B}(\mathsf{X}) \; . 
\end{equation}
In this case, $C$ is said to be \textit{small}.  If \eqref{eq:min_con}
holds with $C=\mathsf{X}$ then $X$ is uniformly ergodic and, as is
well-known,
\[
\|P^{n}(x,\cdot) - \pi(\cdot) \| \le (1-\epsilon)^{\lfloor
  n/n_{0}\rfloor} \; .
\]
Uniformly ergodic Markov chains are rarely encountered in MCMC unless
$\mathsf{X}$ is finite or bounded. 

Geometric ergodicity may be established via the following drift
condition: Suppose that for a function $V : \mathsf{X} \rightarrow
[1,\infty)$ there exist constants $d >0$, $b < \infty$ such that
\begin{equation}
\label{eq:geo_dc}
\Delta V(x) \le -d V(x) + b I(x \in C) \hspace*{5mm} x \in
\mathsf{X} 
\end{equation}
where $C$ is a small set and $I$ is the usual indicator function.

Polynomial ergodicity may be established via a slightly different
drift condition: Suppose that for a function $V : \mathsf{X}
\rightarrow [1,\infty)$ there exist constants $d > 0$, $b<\infty$ and
$0\le \tau <1$ such that
\begin{equation}
\label{eq:poly_dc}
\Delta V(x) \le -d [V(x)]^{\tau} + b I(x \in C) \hspace*{5mm} x \in
\mathsf{X} 
\end{equation}
where $C$ is a small set.  \citet{jarn:robe:2002} show that
\eqref{eq:poly_dc} implies that $X$ is polynomially ergodic of degree
$\tau/(1-\tau)$.  \citet{douc:fort:moul:soul:2004} have recently
generalized this drift condition to other subgeometric (slower than
geometric) rates of convergence.

\begin{remark}
  Either of the drift conditions \eqref{eq:geo_dc} or
  \eqref{eq:poly_dc} imply that in \eqref{eq:tvbd} we can take $M(x)
  \propto V(x)$.  Moreover, Theorem 14.3.7 in \citet{meyn:twee:1993}
  shows that if \eqref{eq:geo_dc} holds then $\text{E}_{\pi} V <
  \infty$.  Since geometric ergodicity is equivalent to
  \eqref{eq:geo_dc} \citep[][Chapter 16]{meyn:twee:1993} we conclude
  that geometrically (and uniformly) ergodic Markov chains satisfy
  \eqref{eq:tvbd} with $\text{E}_{\pi} M < \infty$.  On the other
  hand, the polynomial drift \eqref{eq:poly_dc} only seems to imply
  that $\text{E}_{\pi} V^{\tau} < \infty$ where $\tau < 1$.  Thus,
  when \eqref{eq:poly_dc} holds, to ensure that $\text{E}_{\pi} M <
  \infty$ we will have to show that $\text{E}_{\pi} V < \infty$.
\end{remark}

Beyond establishing a rate of convergence, drift conditions also
immediately imply the existence of a CLT for certain functions.

\begin{theorem}
\label{thm:dm_clt}
Let $X$ be a Harris ergodic Markov chain on $\mathsf{X}$ having
stationary distribution $\pi$.  Suppose $f : \mathsf{X} \rightarrow
\mathbb{R}$ and assume that one of the following conditions hold:
\begin{enumerate}
\item The drift condition \eqref{eq:geo_dc} holds and $f^{2} (x) \le
  V(x)$ for all $x \in \mathsf{X}$.
\item The drift condition \eqref{eq:poly_dc} holds and $|f(x)| \le
  V(x)^{\tau + \eta -1}$ for all $x \in \mathsf{X}$ where $1-\tau \le
  \eta \le 1$ is such that $\text{E}_{\pi} V^{2\eta} < \infty$.
\end{enumerate} 
Then $\sigma_{f}^{2} \in [0, \infty)$ and if $\sigma^{2}_{f} > 0$ then
for any initial distribution
\begin{equation*}
\sqrt{n} (\bar{f}_{n} - \text{E}_{\pi} f) \stackrel{d}{\rightarrow}
\text{N} (0, \sigma^{2}_{f}) 
\end{equation*}
as $n \rightarrow \infty$.
\end{theorem}

\begin{remark}
  The first part of the theorem is from \citet[Theorem
  17.0.1]{meyn:twee:1993} while the second part is due to
  \citet[Theorem 4.2]{jarn:robe:2002}.  
\end{remark}

\begin{remark}
\citet{kont:meyn:2003}
  investigate the rate of convergence in the CLT when the drift
  condition \eqref{eq:geo_dc} holds.
\end{remark}

There has been a substantial amount of effort devoted to establishing
drift and minorization conditions in MCMC settings.  For example,
\citet{hobe:geye:1998}, \citet{jone:hobe:2004},
\citet{marc:hobe:2004}, \citet{robe:1995}, \citet{robe:pols:1994},
\citet{rose:1995a,rose:1996} and \citet{tier:1994} examined Gibbs
samplers while \citet{chri:moll:waag:2001},
\citet{douc:fort:moul:soul:2004}, \citet{fort:moul:2000},
\citet{fort:moul:2003}, \citet{geye:1999}, \citet{jarn:hans:2000},
\citet{jarn:robe:2002}, \citet{meyn:twee:1994b}, and
\citet{meng:twee:1996} considered Metropolis-Hastings-Green (MHG)
algorithms.  Also, \citet{mira:tier:2002} and \citet{robe:rose:1999a}
worked with slice samplers.

In the next section three simple examples are presented in order to give
the reader a taste of using these results in specific models and to
demonstrate the application of Theorem~\ref{thm:dm_clt}.  More
substantial examples will be considered in Section~\ref{sec:ex}.

\subsection{Examples}
\label{sec:toy_ex}

\noindent {\textit{Example} 1 continued.}
Since $\mathsf{X}$ is finite it is easy to see that \eqref{eq:min_con}
holds with $C=\mathsf{X}$ and hence the Markov chain described in
Example~\ref{ex:hardcore} is uniformly ergodic.  Of course, if $n_{1}$
and $n_{2}$ are reasonably large $\epsilon$ may be too small to be
useful. \medskip

\begin{example}
\label{ex:dmc}
Suppose $X$ lives on $\mathsf{X}=\mathbb{Z}$ such that if $x \ge 1$
and $0 < \theta < 1$ then
\[
P(x, x+1)  = P(-x, -x-1) = \theta \, ,  \hspace{5mm} 
P(x,0)  =  P(-x, 0) =1 - \theta\, ,  
\]
\[
P(0,1)  = P(0, -1) = \frac{1}{2} \; .
\]
This chain is Harris ergodic and has stationary distribution given by
$\pi(0) = (1-\theta)/(2-\theta)$ and for $x \ge 1$
\[
\pi(x)=\pi(-x)=\pi(0) \frac{\theta^{x-1}}{2} \; .
\]
In Appendix~\ref{app:dmc} the drift condition \eqref{eq:geo_dc} is
verified with $V(x) =a^{|x|}$ for $a > 1$ satisfying $a \theta < 1$
and $(a \theta -1)a + 1 - \theta < 0$ and $C=\{0\}$.  Hence a CLT
holds for $\bar{f}_{n}$ if $f^{2}(x) \le a^{|x|} \; \forall \, x \in
\mathbb{Z}$.
\end{example}

\begin{example}
\label{ex:random walk}
\citet{jarn:robe:2002} and \citet{tuom:twee:1994} consider the
following example and establish a polynomial rate of convergence.  Let
$X$ be a random walk on $[0,\infty)$ determined by
\[
X_{n+1} = (X_{n} + W_{n+1})^{+} 
\]
where $W_{1}, W_{2}, \ldots$ is a sequence of independent and
identically distributed real-valued random variables.  As long as
$\text{E} (W_{1}) < 0$ this chain will be Harris ergodic.  When
$\text{E} (W_{1}^{+})^{m} < \infty$ for some $m \ge 2$
\citet{jarn:robe:2002} establish the drift condition
\eqref{eq:poly_dc} with $V(x) = (x+1)^{m}$, $\tau= (m-1)/m$ and
$C=[0,k]$ for some $k < \infty$. Hence a CLT holds for $\bar{f}_{n}$
if $|f(x)| \le (x+1)^{m(\tau+\eta-1)}$ for all $x \ge 0$ where $1-\tau
\le \eta \le 1$ is such that $\text{E}_{\pi} (x+1)^{2m\eta} < \infty$.
Note that this moment condition also implies that $\text{E}_{\pi} V <
\infty$ as long as $\eta \ge 1/2$.  Hence by an earlier remark
$\text{E}_{\pi} M < \infty$ with $M$ as in \eqref{eq:tvbd}.
\end{example}\bigskip

Two things are clear: (i) drift and minorization provide powerful
constructive tools for establishing a rate of convergence in total
variation; and (ii) they are less impressive (but often useful!) tools
for establishing CLTs in that the results in Theorem~\ref{thm:dm_clt}
depend on the non-unique function $V$.

\section{Mixing Sequences}
\label{sec:mixing}
The goal of this section is to introduce three types of mixing
conditions and discuss some of the connections with the total
variation convergence in \eqref{eq:bas_conv} and \eqref{eq:tvbd}.
There are a variety of mixing conditions (e.g. absolute regularity)
that will not be considered here since they don't seem to have much
impact on the CLT.  Roughly speaking, mixing conditions are all
attempts to quantify the rate at which events in the distant future
become independent of the past.

Let $Y:=\{Y_{n}\}$ denote a general sequence of random variables on a
probability space $(\Omega, {\cal F}, {\cal P})$ and let ${\cal
  F}_{k}^{m} = \sigma(Y_k, \ldots, Y_m)$. 
  
\begin{definition}
  The sequence $Y$ is said to be strongly mixing (or $\alpha$--mixing)
  if $\alpha(n) \rightarrow 0$ as $n \rightarrow \infty$ where
\[
\alpha(n) := \sup_{k \ge 1} \sup_{A \in {\cal F}_{1}^{k}, B \in {\cal
    F}_{k+n}^{\infty}} | {\cal P}(A \cap B) - {\cal P}(A)
    {\cal P}(B) | \; .  
\]
\end{definition}

Harris ergodic Markov chains are strongly mixing.  Recall the coupling
inequality \citep[p.  12][]{lind:1992}:
\begin{equation}
\label{eq:ci}
\| P^{n} (x, \cdot) - \pi(\cdot) \| \le \text{Pr}_{x} (T>n) 
\end{equation}
where $T$ is the usual coupling time of two Markov chains; one started
in stationarity and one started arbitrarily.  Under our assumptions
the coupling time is almost surely finite and $\Pr(T>n) \rightarrow 0$
as $n \rightarrow \infty$. Let $A$ and $B$ be Borel sets so that by
\eqref{eq:ci}
\[
| P^{n} (x, A) - \pi(A) | \le  \text{Pr}_{x}(T > n)
\]
and
\[
\begin{split}
\int_{B} \text{Pr}_{x} (T > n) \pi(dx) & \ge \int_{B} | P^{n} (x, A) -
\pi(A) | \pi(dx) \\  
& \ge | \int_{B} [P^{n} (x, A) - \pi(A)]  \pi(dx) |\\
& = | \Pr (X_{n} \in A \text{ and } X_{0} \in B) - \pi(A) \pi(B) | \; .
\end{split}
\]
Then $\alpha(n) \le \text{E}_{\pi}[ \text{Pr}_{x} (T > n)]$ and a
dominated convergence argument shows that $\text{E}_{\pi}[
\text{Pr}_{x} (T > n)] \rightarrow 0$ as $n \rightarrow \infty$ and
hence $\alpha(n) \rightarrow 0$ as $n \rightarrow \infty$.  Moreover,
the rate of total variation convergence bounds the rate of
$\alpha$--mixing: if \eqref{eq:tvbd} holds with $\text{E}_{\pi} M <
\infty$, a similar argument shows that $\alpha(n) \le \gamma(n)
\text{E}_{\pi} M$ and hence $\alpha(n) = O(\gamma(n))$.  For example,
geometrically ergodic Markov chains enjoy exponentially fast strong
mixing.

Suppose the process $Y$ is strictly stationary and let $f : \mathbb{R}
\rightarrow \mathbb{R}$ be a Borel function.  Define the process $W :=
\{W_n = f(Y_n)\}$.  Set ${\cal G}_{k}^{m} := \sigma(W_k , \ldots,
W_m)$; hence ${\cal G}_{k}^{m} \subseteq {\cal F}_{k}^{m}$.  Let
$\alpha_{W}$ and $\alpha_{Y}$ be the strong mixing coefficients for
the processes $W$ and $Y$, respectively.  Then $\alpha_{W} (n) \le
\alpha_{Y} (n)$.  Similar comments apply to the mixing conditions
given below. This elementary observation is fundamental to the proofs
of the Markov chain CLTs considered in the sequel.

\begin{definition}
  The sequence $Y$ is said to be asymptotically uncorrelated (or
  $\rho$--mixing) if $\rho(n) \rightarrow 0$ as $n \rightarrow \infty$
  where
\[
\rho(n) := \sup \{\text{corr}(U,V), \, U \in L_{2}({\cal F}_{1}^{k})
\, , V \in L_{2}({\cal F}_{k+n}^{\infty}) \, k \ge 1\} .
\]
\end{definition}
It is standard that $\rho$-mixing sequences are also strongly mixing
and, in fact, $4 \alpha(n) \le \rho(n)$.  It is a consequence of the
strong Markov property that if a Harris ergodic Markov chain is
$\rho$-mixing then it enjoys exponentially fast $\rho$-mixing
\citep[Theorem 4.2]{brad:1986} in the sense that there exists a
$\theta >0$ such that $\rho(n) = O(e^{-\theta n})$.  

\citet{rose:1971} develops a necessary and sufficient condition for a
Markov chain to be $\rho$--mixing but before giving it a slight
digression is required.  Define the Hilbert space $L^{2}(\pi) := \{
f\, : \mathsf{X} \rightarrow \mathbb{R} \, ; \, \text{E}_{\pi} f^{2} <
\infty \}$ with inner product $(f,g) = \text{E}_{\pi} [f(x) g(x)] $
and norm $\|\cdot\|_{2}$.  Let $L^{2}_{0}(\pi) := \{f \in L^{2}(\pi)
\, ;\, \text{E}_{\pi} f =0\}$ and note that if $f, g \in
L^{2}_{0}(\pi)$ then $(f,g) = \text{cov}_{\pi} (f,g)$.  The kernel $P$
defines an operator $T : L^{2}(\pi) \rightarrow L^{2}(\pi)$ via
\[
(Tf)(x) = \int P(x, dy) f(y) \; .
\]
It is easy to show that $T$ is a contraction (i.e., $\|T\| \le 1$).
Also, $T$ is self-adjoint if and only if the kernel $P$ satisfies 
\textit{detailed balance} with respect to $\pi$:
\begin{equation}
\label{eq:dbc}
\pi(dx) P(x,dy) = \pi(dy) P(y, dx) \hspace{5mm} \forall \, x,y \in
\mathsf{X} \; .
\end{equation}
\citet[p. 207]{rose:1971} shows that a Harris ergodic Markov chain is
$\rho$--mixing if and only if 
\begin{equation}
\label{eq:rhoeq}
\lim_{n \rightarrow \infty} \sup_{\stackrel{f \in
  L^{2}_{0}(\pi)}{\|f\|_{2}=1}} \|T^{n} f \|_{2} = 0 \; .
\end{equation}
There has been some work done on establishing sufficient conditions
for Markov chains to be $\rho$-mixing.  For example,
\citet{liu:wong:kong:1995} show that if the operator induced by a
Gibbs sampler satisfies a Hilbert--Schmidt condition then it is
$\rho$-mixing.  However, the most interesting case is given by
\citet{robe:rose:1997c} whose Theorem 2.1 shows that if $X$ is
geometrically ergodic and \eqref{eq:dbc} holds then there exists a $c
< 1$ such that $\| T f \|_{2} \le c^{2}$ and $\| T^{n} f \|_{2} = \| T
f \|_{2}^{n}$ hence \eqref{eq:rhoeq} holds.  We conclude that if $X$
is geometrically ergodic and \eqref{eq:dbc} holds then $X$ is
asymptotically uncorrelated.

\begin{remark}
  Many Markov chains satisfy \eqref{eq:dbc}, indeed the MHG algorithm
  satisfies \eqref{eq:dbc} by construction.  However, \eqref{eq:dbc}
  does not hold for those Markov chains associated with systematic
  scan Gibbs samplers and the Markov chain in Example~\ref{ex:dmc},
  for example.
\end{remark} 

\begin{definition}
  The sequence $Y$ is said to be uniformly mixing (or $\phi$--mixing)
  if $\phi(n) \rightarrow 0$ as $n \rightarrow \infty$ where
\[
\phi(n) := \sup_{k \ge 1} \sup_{\stackrel{A \in {\cal F}_{1}^{k},
    {\cal P}(A) \neq 0}{B \in {\cal  F}_{k+n}^{\infty}}} | {\cal P}(B
    | A) - {\cal P}(B) | \; . 
\]
\end{definition}
Uniformly mixing sequences are also asymptotically uncorrelated and
strongly mixing.  Moreover, $\rho(n) \le 2 \sqrt{\phi(n)}$.  A Harris
ergodic Markov chain is uniformly ergodic if and only if it is
uniformly mixing; see \citet[][pp 367--368]{ibra:linn:1971}.

As with asymptotically uncorrelated sequences it is a consequence of
the strong Markov property that if a Harris ergodic Markov chain is
$\phi$-mixing then it enjoys exponentially fast $\phi$-mixing
\citep[Theorem 4.2]{brad:1986} in the sense that there exists a
$\theta >0$ such that $\phi(n) = O(e^{-\theta n})$.

We collect and concisely state the main conclusions of this section.

\begin{theorem}
Let $X$ be a Harris ergodic Markov chain with stationary distribution
$\pi$. 
\begin{enumerate}
\item $X$ is strongly mixing, i.e., $\alpha(n) \rightarrow 0$.
\item If \eqref{eq:tvbd} holds with $\text{E}_{\pi} M < \infty$ then
  $\alpha(n) = O(\gamma(n))$.
\item If $X$ is geometrically ergodic and \eqref{eq:dbc} holds then $X$
  is asymptotically uncorrelated, in which case there exists a
  $\theta>0$ such that $\rho(n) = O(e^{-\theta n})$.
\item $X$ is uniformly ergodic if and only if $X$ is uniformly mixing,
  in which case there exists a $\theta>0$ such that $\phi(n) =
  O(e^{-\theta n})$. 
\end{enumerate}
\end{theorem}

\section{Central Limit Theorems}
\label{sec:clt}
We begin with a characterization of the CLT for strongly mixing
processes.  Define $S_n = \sum_{i=1}^{n} Y_i$ and $\sigma_{n}^{2} = E
S_{n}^{2}$.

\begin{theorem} \label{thm:t1} \citep{cogb:1960,denk:1986,mori:yosh:1986}
  Let $Y$ be a centered strictly stationary strongly mixing sequence
  such that $E Y_{0}^{2} < \infty$.  If $\sigma_{n}^{2} \rightarrow
  \infty$ as $n \rightarrow \infty$ then the following are equivalent:
\begin{enumerate}
\item $S_{n} / \sigma_{n} \stackrel{d}{\rightarrow} \text{N}(0,1)$.
\item $\{ S_{n}^{2} / \sigma_{n}^{2} \, , \, n \ge 1 \}$ is uniformly
  integrable. 
\end{enumerate}
\end{theorem}\bigskip

\begin{remark}
  Since Harris ergodic Markov chains are strongly mixing this result
  is applicable in MCMC settings.
\end{remark}

\begin{remark}  \label{rem:init}
The assumption of stationarity is not an issue for Harris ergodic
Markov chains since if a CLT holds for any one initial distribution
then it holds for every initial distribution \citep[Proposition
17.1.6]{meyn:twee:1993}.  
\end{remark}\medskip

\noindent \citet{chen:1999} provides the following characterization of
the CLT. 

\begin{theorem} \label{thm:t2} \citep{chen:1999}
  Let $X$ be a Harris ergodic Markov chain and $f$ be a function such
  that $\text{E}_{\pi} f =0$ and $\text{E}_{\pi} f^{2} < \infty$.
  Then the following are equivalent:
\begin{enumerate}
\item $\sqrt{n} \bar{f}_{n} \stackrel{d}{\rightarrow}
  \text{N}(0,\sigma^{2})$ for some $\sigma^{2} \ge 0$.
\item $\{\sqrt{n} \bar{f}_{n}\, , \, n \ge 1\}$ is bounded in probability.
\end{enumerate}
\end{theorem}\bigskip

\begin{remark}
\citet{chen:1999} also provides another equivalent condition in terms
of quantities based on the so-called \textit{split chain}.  But this
is not germane to the current discussion. 
\end{remark}

\subsection{Sufficient Conditions}
\label{sec:suff} 

\begin{theorem} \citep{ibra:1962,ibra:linn:1971} \label{thm:il} 
  Let $Y$ be a centered strictly stationary strongly mixing sequence.
  Suppose at least one of the following conditions:
\begin{enumerate}
\item There exists $B < \infty$ such that $|Y_{n}| < B \; a.s.$ and
  $\sum_{n} \alpha(n) < \infty \; ; \; \text{or}$ 
\item $\text{E} |Y_n|^{2 + \delta} < \infty$ for some $\delta >0$ and
\begin{equation}
\label{eq:alpha}
\sum_{n} \alpha(n)^{\delta/(2+\delta)} < \infty \; .
\end{equation}
\end{enumerate}
Then
\[
\sigma^{2} = \text{E}(Y_{0}^{2}) + 2 \sum_{j=1}^{\infty} \text{E}
(Y_{0} Y_{j}) < \infty
\]
and if $\sigma^{2} >0$, as $n \rightarrow \infty$,
\[
n^{-1/2} S_n \stackrel{d}{\rightarrow} \text{N} (0, \sigma^{2}) \; .
\]
\end{theorem}

\begin{corollary} \label{cor:tv_clt}
  Let $f : \mathsf{X} \rightarrow \mathbb{R}$ be a Borel function such
  that $\text{E}_{\pi} |f(x)|^{2 + \delta} < \infty$ for some $\delta
  > 0$ and suppose $X$ is a Harris ergodic Markov chain with
  stationary distribution $\pi$.  If \eqref{eq:tvbd} holds such that
  $\text{E}_{\pi} M < \infty$ and $\gamma(n)$ satisfies
\begin{equation}
\label{eq:gamma_sum}
\sum_{n} \gamma(n)^{\delta/(2+\delta)} < \infty 
\end{equation}
then for any initial distribution, as $n \rightarrow \infty$
\[
\sqrt{n} (\bar{f}_{n} -  \text{E}_{\pi} f) \stackrel{d}{\rightarrow}
\text{N} (0, \sigma_{f}^{2}) \; .
\]
\end{corollary}
Later, CLTs for $\phi$-mixing and $\rho$-mixing Markov chains will be
presented.  However, the proofs of these results are similar to the
proof of Corollary~\ref{cor:tv_clt}.  Hence only the following proof
is included.
\begin{proof}
  Let $\alpha(n)$ and $\alpha_{f}(n)$ denote the strong mixing
  coefficients for the Markov chain $X=\{ X_{n} \}$ and the functional
  process $\{ f(X_{n}) \}$, respectively.  By an earlier remark
  $\alpha_{f} (n) \le \alpha(n)$ for all $n \ge 1$.  Moreover, we have
  that $\alpha(n) \le \gamma(n) \text{E}_{\pi} M$ where $\gamma(n)$
  and $M$ are given in \eqref{eq:tvbd}.  Hence \eqref{eq:gamma_sum}
  guarantees that
\[ 
\sum_{n} \alpha_{f}(n)^{\delta/(2+\delta)} < \infty 
\]
and the result follows from the Theorem and Remark~\ref{rem:init}.
\end{proof}

Corollary~\ref{cor:tv_clt} immediately yields some special cases which
have proven to be useful in MCMC settings.

\begin{corollary} \label{cor:ge_pe}
  Suppose $X$ is a Harris ergodic Markov chain with stationary
  distribution $\pi$ and let $f : \mathsf{X} \rightarrow \mathbb{R}$ be a
  Borel function. Assume one of the following conditions:
\begin{enumerate}
\item \citep{chan:geye:1994} $X$ is geometrically ergodic and
  $\text{E}_{\pi} |f(x)|^{2+\delta} < \infty$ for some $\delta > 0$;
\item $X$ is polynomially ergodic of order $m$, $\text{E}_{\pi} M <
  \infty$ and $\text{E}_{\pi} |f(x)|^{2+\delta} < \infty$ where $m\delta
  > 2+\delta$; or
\item $X$ is polynomially ergodic of order $m > 1$, $\text{E}_{\pi} M
  < \infty$ and there exists $B< \infty$ such that $|f(x)| < B$
  $\pi$-almost surely.
\end{enumerate}
Then for any initial distribution, as $n \rightarrow \infty$
\[
\sqrt{n} (\bar{f}_{n} - \text{E}_{\pi} f) \stackrel{d}{\rightarrow}
\text{N} (0, \sigma_{f}^{2}) \; .
\]
\end{corollary}

For geometrically ergodic Markov chains the moment condition can not
be weakened to a second moment (i.e., $\text{E}_{\pi} f^{2}(x) <
\infty$) without additional assumptions.  \citet{hagg:2004} has
recently established the existence of a geometrically ergodic Markov
chain and a function $f$ such that $\text{E}_{\pi} f^{2}(x) < \infty$
yet a CLT fails for any choice of $\sigma^{2}$.  Also, see
\citet{brad:1983} for a counterexample with the same conclusion.
These results are not too surprising since there are non-trivial
counterexamples that indicate that the conditions of
Theorem~\ref{thm:il} are nearly as good as can be expected.  For
example, \citet{hern:1983} constructs a strictly stationary sequence
of uncorrelated random variables, $\{ Y_{n} \}$, that have an
arbitrarily fast strong mixing rate and $0 < E Y_{1}^{2} < \infty$ yet
the CLT fails.  Further counterexamples have been given by
\citet{davy:1973} and \citet{brad:1985}.  However, a slightly weaker
moment condition is available if the sequence enjoys at least
exponentially fast strong mixing which is the case for geometrically
ergodic Markov chains.  The following theorem is a special case of a
result in \citet{douk:mass:rio:1994}.

\begin{theorem}
\citep{douk:mass:rio:1994} \label{thm:dmr} 
  Let $Y$ be a centered strictly stationary strongly mixing sequence.
  If the strong mixing coefficients satisfy $\alpha(n) = O(a^{n})$ for
  some $ 0 < a < 1$ and $\text{E} [Y_{0}^{2} (\log^{+} |Y_{0}|) <
  \infty$ then 
\[
\sigma^{2} = \text{E} Y_{0}^{2} + 2 \sum_{k=1}^{\infty} \text{E} (Y_{0}
Y_{k})
\]
converges absolutely and if $\sigma^{2} > 0$, as $n \rightarrow
\infty$ 
\[
n^{-1/2} S_n \stackrel{d}{\rightarrow} \text{N} (0,\sigma^{2})
\; .
\] 
\end{theorem}

\begin{corollary} \label{cor:ge}
  Suppose $X$ is a Harris ergodic Markov chain with stationary
  distribution $\pi$ and let $f : \mathsf{X} \rightarrow \mathbb{R}$
  be a Borel function.  If $X$ is geometrically ergodic and
  $\text{E}_{\pi} [f^{2} (x) (\log^{+} |f(x)|)] < \infty$ then for any
  initial distribution, as $n \rightarrow \infty$
\[
\sqrt{n} (\bar{f}_{n} - \text{E}_{\pi} f) \stackrel{d}{\rightarrow}
\text{N} (0, \sigma_{f}^{2}) \; .
\]
\end{corollary}

A weaker moment condition is available for $\rho$--mixing sequences. 

\begin{theorem} \citep{ibra:1975} \label{thm:pell}
  Let $Y$ be a centered strictly stationary $\rho$--mixing sequence
  with $\text{E}\, Y_{0}^{2} < \infty$.  Suppose
\begin{equation}
\label{eq:rho_rate}
\sum_{n=1}^{\infty} \rho(n) < \infty \; .
\end{equation}
Then 
\[
\sigma^{2} = \text{E} Y_{0}^{2} + 2 \sum_{k=1}^{\infty} \text{E} (Y_{0}
Y_{k})
\]
converges absolutely and if $\sigma^{2} > 0$, as $n \rightarrow
\infty$ 
\[
n^{-1/2} S_n \stackrel{d}{\rightarrow} \text{N} (0,\sigma^{2})
\; .
\]
\end{theorem}

\noindent Recall that if the Markov chain $X$ is geometrically ergodic
and satisfies detailed balance, it enjoys exponentially fast
$\rho$--mixing and hence \eqref{eq:rho_rate} obtains.

\begin{corollary} \citep{robe:rose:1997c} \label{cor:rr97}  
  Let $X$ be a geometrically ergodic Markov chain with stationary
  distribution $\pi$.  Suppose $X$ satisfies \eqref{eq:dbc} and that
  $\text{E}_{\pi} f^{2}(x) < \infty$.  Then for any initial distribution,
  as $n \rightarrow \infty$
\[
\sqrt{n} (\bar{f}_{n} -  \text{E}_{\pi} f) \stackrel{d}{\rightarrow}
\text{N} (0, \sigma_{f}^{2}) \; .
\]
\end{corollary}

\begin{remark}
  \citet{robe:rose:1997c} obtain this result via Corollary 1.5 in
  \citet{kipn:vara:1986}.  We have thus provided an alternative
  derivation.
\end{remark}

An accessible proof of the following result may be found in
\citet{bill:1968} and \citet{ibra:linn:1971}.  Also see Chapter 5 of
\citet{doob:1953} and Lemma 3.3 in \citet{cogb:1972}.

\begin{theorem} \label{thm:bill68}
  Let $Y$ be a centered strictly stationary uniformly mixing sequence
  with $\text{E}\, Y_{0}^{2} < \infty$.  If
\begin{equation}
\label{eq:phi_rate}
\sum_{n} \sqrt{\phi(n)} < \infty
\end{equation}
then
\[
\sigma^{2} = \text{E} Y_{0}^{2} + 2 \sum_{k=1}^{\infty} \text{E} (Y_{0}
Y_{k})
\]
converges absolutely and if $\sigma^{2} > 0$ then as $n \rightarrow
\infty$
\[
n^{-1/2} S_n \stackrel{d}{\rightarrow} \text{N} (0,\sigma^{2})
\; .
\]
\end{theorem}

\noindent If $X$ is uniformly ergodic the coefficients $\phi(n)$
decrease exponentially and \eqref{eq:phi_rate} is obvious.

\begin{corollary} \label{cor:t94} \citep{ibra:linn:1971,tier:1994} 
  Let $X$ be a uniformly ergodic Markov chain with stationary
  distribution $\pi$.  Suppose $\text{E}_{\pi} f^{2}(x) < \infty$.
  Then for any initial distribution, as $n \rightarrow \infty$
\[
\sqrt{n} (\bar{f}_{n} -  \text{E}_{\pi} f) \stackrel{d}{\rightarrow}
\text{N} (0, \sigma_{f}^{2}) \; .
\]
\end{corollary}

The main conclusions of this section can be concisely stated as
follows. 
\begin{theorem}
Let $X$ be a Harris ergodic Markov chain on $\mathsf{X}$ with invariant
distribution $\pi$ and let  $f : \mathsf{X} \rightarrow \mathbb{R}$ is
a Borel function. Assume one of the following conditions:
\begin{enumerate}
\item $X$ is polynomially ergodic of order $m > 1$, $\text{E}_{\pi} M
  < \infty$ and there exists $B< \infty$ such that $|f(x)| < B$ almost
  surely;
\item $X$ is polynomially ergodic of order $m$, $\text{E}_{\pi} M <
  \infty$ and $\text{E}_{\pi} |f(x)|^{2+\delta} < \infty$ where
  $m\delta > 2+\delta$;
\item $X$ is geometrically ergodic and $\text{E}_{\pi}
  |f(x)|^{2+\delta} < \infty$ for some $\delta > 0$;
\item $X$ is geometrically ergodic and $\text{E}_{\pi} [f^{2}(x)
  (\log^{+} |f(x)|)] 
  < \infty$; 
\item $X$ is geometrically ergodic, satisfies \eqref{eq:dbc} and
  $\text{E}_{\pi} f^{2}(x) < \infty$; or
\item $X$ is uniformly ergodic and $\text{E}_{\pi} f^{2}(x) < \infty$.
\end{enumerate}
Then for any initial distribution, as $n \rightarrow \infty$
\[
\sqrt{n} (\bar{f}_{n} - \text{E}_{\pi} f) \stackrel{d}{\rightarrow}
\text{N} (0, \sigma_{f}^{2}) \; .
\]
\end{theorem}

\begin{remark}
  Condition 1 of the theorem is interesting for applications of MCMC
  in Bayesian settings.  In this case, it is often the case that
  posterior probabilities, i.e. expectations of indicator functions,
  are of interest.  Since indicator functions are bounded it follows
  that a CLT will hold under a very weak mixing condition.
\end{remark}

\section{Examples}
\label{sec:ex}

\subsection{Toy Examples Revisited}
\label{sec:toy}

\noindent {\textit{Example} 1 continued.}  Recall that since
$\mathsf{X}$ is finite this chain is uniformly ergodic and uniformly
mixing.  Hence Corollary~\ref{cor:t94} implies that a CLT will hold
for $\bar{f}_{n}$ if $\text{E}_{\pi} f^{2}(x) < \infty$ which will
hold except in unusual cases.\smallskip

\noindent {\textit{Example} 2 continued.}  This chain is geometrically
ergodic but does not satisfy \eqref{eq:dbc}.  Hence it is strongly
mixing and we can not conclude that it is asymptotically uncorrelated.
Thus the best we can do is to appeal to Corollary~\ref{cor:ge} and
conclude that a CLT will hold for $\bar{f}_{n}$ if $\text{E}_{\pi}
[f(x)^{2} (\log^{+} |f(x)|)] < \infty$.  Recall that in
subsection~\ref{sec:toy_ex} it was shown that a CLT holds for
$\bar{f}_{n}$ if $f^{2}(x) \le a^{|x|} \; \forall \, x \in \mathbb{Z}$
when $a > 1$ satisfies $a \theta < 1$ and $(a \theta -1)a + 1 - \theta
< 0$.\smallskip

\noindent {\textit{Example} 3 continued.}  Let $m > 2$ and recall that
this random walk is polynomially ergodic of order $m-1$ and that
Theorem~\ref{thm:dm_clt} says a CLT holds if $f(x) \le
(x+1)^{m(\tau+\eta-1)}$ for all $x \ge 0$ where $1-\tau \le \eta \le
1$ is such that $\text{E}_{\pi} (x+1)^{2m\eta} < \infty$.
Alternatively, an appeal to Corollary~\ref{cor:ge_pe} says that we
have a CLT if $\text{E}_{\pi} (x+1)^{m} < \infty$ and $\text{E}_{\pi}
|f(x)|^{2+\delta} < \infty$ where $\delta > 2/(m-2)$.

\subsection{A Benchmark Gibbs Sampler}
\label{sec:gamma}
The following Gibbs sampler is similar to one used by many authors to
analyze the benchmark pump failure data set in \citet{gave:o'mu:1987}.
For example, \citet{robe:case:1999}, \citet{rose:1995a}, and
\citet{tier:1994} consider similar settings and establish uniform
ergodicity of the corresponding Gibbs samplers.  

Set $y=(y_{1} , y_{2}, \ldots, y_{n})^{T}$ and let $\pi(x,y)$ be a
joint density on $\mathbb{R}^{n+1}$ such that the corresponding full
conditionals are
\[
\begin{split}
X | y & \sim \text{Gamma} \, (\alpha_{1}, a+b^{T} y) \\
Y_{i} | x & \sim \text{Gamma} \, (\alpha_{2i} , \beta_{i}(x))
\end{split}
\]
for $i=1,\ldots,n$, $b=(b_{1},\ldots,b_{n})^{T}$ where $a > 0$ and
each $b_{i} > 0$ are known.  (Say $U \sim \text{Gamma}(\alpha, \beta)$
if its density is proportional to $u^{\alpha-1} e^{-u \beta} I(u >
0)$.)  Since, conditional on $x$, the order in which the $Y_{i}$ are
updated is irrelevant we will use a two variable Gibbs sampler with
the transition rule $(x', y') \rightarrow (x,y)$; that is, given that
the current value is $(X_{n}=x', Y_{n}=y')$ we obtain $(X_{n+1} ,
Y_{n+1})$ by first drawing $x \sim f(X_{n+1} | y')$ then $Y_{i, n+1}
\sim f(y_{i} | x)$.  A minor modification of Tierney's argument will
show that \eqref{eq:min_con} holds on $C=\mathsf{X}$ with $n_{0}=1$
and if for $i=1,\ldots,n$ there is a function $g >0$ such that for all
$x > 0$
\begin{equation}
\label{eq:gamma_ineq}
\frac{\beta_{i}(x)}{b_{i} x + \beta_{i}(x)} \ge g(x) \; .
\end{equation}
Thus if \eqref{eq:gamma_ineq} holds this Gibbs sampler is uniformly
ergodic (or uniformly mixing) and an appeal to Corollary~\ref{cor:t94}
shows that a CLT is assured if $\text{E}_{\pi} f^{2}(x) < \infty$.

\subsection{A Gibbs Sampler for a Hierarchical Bayes Setting}
\label{sec:bayes}
Consider the following Bayesian version of the classical normal theory
one--way random effects model.  First, conditional on
$\theta=(\theta_1,\dots,\theta_K)^T$ and $\lambda_e$, the data,
$Y_{ij}$, are independent with
\[
Y_{ij}|\theta,\lambda_{e} \sim
\mbox{N}(\theta_{i},\lambda_{e}^{-1})
\]
where $i=1,\ldots,K$ and $j=1,\ldots,m_{i}$.  Conditional on $\mu$ and
$\lambda_{\theta}$, $\theta_1, \ldots, \theta_K$ and $\lambda_e$ are
independent with
\[
\theta_i | \mu,\lambda_{\theta} \sim \mbox{N}(\mu ,
 \lambda_{\theta}^{-1} ) \hspace*{5mm}\mbox{and} \hspace*{5mm}
 \lambda_{e} \sim \mbox{Gamma}(a_{2},b_{2}),
\]
where $a_2$ and $b_2$ are known positive constants.  Finally, $\mu$
and $\lambda_{\theta}$ are assumed independent with 
\[
\mu \sim \mbox{N}(m_0,s_0^{-1}) \hspace*{5mm}\mbox{and}
 \hspace*{5mm} \lambda_{\theta} \sim  \mbox{Gamma}(a_{1},b_{1})
\]
where $m_0,s_{0},a_{1}$ and $b_{1}$ are known constants; all of the
priors are proper since $s_{0},a_{1}$ and $b_{1}$ are assumed to be
positive and $m_{0} \in \mathbb{R}$.  The posterior density of this
hierarchical model is characterized by
\begin{equation}
\label{eq:post}
\pi(\theta,\mu,\lambda| y) \propto g(y|\theta,\lambda_{e})
g(\theta|\mu,\lambda_{\theta}) g(\lambda_{e}) g(\mu) g(\lambda_{\theta})
\end{equation}
where $\lambda=(\lambda_{\theta},\lambda_{e})^T$, $y$ is a vector
containing all of the data, and $g$ denotes a generic density.
Expectations with respect to $\pi$ typically require evaluation of
ratios of intractable integrals, which may have dimension $K+3$ and
typically, $K \ge 3$.

We are interested in the standard Gibbs sampler which leaves the
posterior \eqref{eq:post} invariant.  Define 
\[
v_1(\theta, \mu) = \sum_{i=1}^K (\theta_{i} - \mu)^{2} \; , \hspace*{7mm}
v_2(\theta) = \sum_{i=1}^K m_{i} ( \theta_{i} -
\overline{y}_{i} )^{2} \hspace*{5mm} \mbox{and} \hspace*{5mm}
\mbox{SSE}=\sum_{i,j}(y_{ij}-\overline{y}_i)^2 \;
\]
where ${\overline y}_{i} = m_i^{-1} \sum_{j=1}^{m_i} y_{ij}$.  The
full conditionals for the variance components are
\begin{equation}
\label{eq:fc_lamt}
\lambda_{\theta} | \theta, \mu,  \lambda_{e}, y \sim
\mbox{Gamma} \left ( \frac{K}{2}+a_1, \frac{v_1(\theta, \mu)}{2} + b_1
\right ) 
\end{equation}
and
\begin{equation}
\label{eq:fc_lame}
\lambda_{e} | \theta, \mu, \lambda_{\theta}, y  \sim \mbox{Gamma} \left (
  \frac{M}{2}+a_2, \frac{v_2(\theta) + \mbox{SSE}}{2} + b_2 \right )
\end{equation}
where $M=\sum_{i} m_{i}$.  If $\theta_{-i} = (\theta_1, \ldots,
\theta_{i-1}, \theta_{i+1}, \dots, \theta_K)^T$ and
$\overline{\theta}=K^{-1} \sum_i \theta_i$, the remaining full
conditionals are
\[
\theta_i | \theta_{-i}, \mu, \lambda_{\theta}, \lambda_{e}, y  \sim
\mbox{ N}\left( \frac{\lambda_{\theta} \mu + m_i \lambda_{e}
    \overline{y}i}{\lambda_{\theta} + m_i \lambda_{e}},
  \frac{1}{\lambda_{\theta} + m_i \lambda_{e}} \right) 
\]
for $i=1,\ldots, K$ and
\[
\mu | \theta, \lambda_{\theta}, \lambda_{e}, y \sim \mbox{ N}\left(
  \frac{s_{0} m_0 + K \lambda_{\theta} \overline{\theta}}{s_{0} +
  K\lambda_{\theta}}, \frac{1}{s_{0} + K \lambda_{\theta}} \right) \,.
\]
Our fixed-scan Gibbs sampler updates $\mu$ then the $\theta_i$'s then
$\lambda_{\theta}$ and $\lambda_{e}$.  Since the $\theta_i$'s are
conditionally independent given $(\mu, \lambda)$, the order in which
they are updated is irrelevant.  The same is true of
$\lambda_{\theta}$ and $\lambda_{e}$ since these two random variables
are conditionally independent given $(\theta,\mu)$.  A one--step
transition of this Gibbs sampler is $(\mu',\theta',\lambda')
\rightarrow (\mu, \theta, \lambda)$ meaning that we sequentially draw
from the conditional distributions $\mu | \lambda' , \theta'$ then
$\theta_{i} | \theta_{-i}, \mu , \lambda'$ for $i=1,\ldots, K$ then
$\lambda_{e} | \mu, \theta$ then $\lambda_{\theta} | \mu, \theta$.
Assume that $m'= \min\{m_1, m_2,\ldots, m_K\} \ge 2$ and that $K \ge
3$.  Let $m''= \max\{m_1,m_2,\ldots, m_K\}$.  Define $\delta_1 = (2
a_1 + K - 2)^{-1}$ and $\delta_2 = (2a_1-2)^{-1} $.

\begin{proposition}
\label{pr:jh_gibbs}
\citep{jone:hobe:2004} Assume that $a_1>3/2$ and $5m'>m''$.  Fix $c_1
\in \big(0, \min\{b_1 , b_2\} \big)$.  Then the Gibbs sampler
satisfies \eqref{eq:geo_dc} with the drift function
\[
V(\theta, \lambda) =  1 + e^{c_1\lambda_{\theta}} + e^{c_1\lambda_{e}} +
  \frac{\delta_2}{K \delta_1 \lambda_{\theta}} +
  \frac{K\lambda_{\theta}}{s_{0} +
  K\lambda_{\theta}}(\overline{\theta} -\overline{y})^2 \; . 
\]
\end{proposition}

\begin{remark}
  \citet{jone:hobe:2004} give values for the constants in
  \eqref{eq:geo_dc} but in an effort to keep the notation under
  control we do not report them here.
\end{remark}
Theorem~\ref{thm:dm_clt} immediately implies a CLT for $\bar{f}_{n}$
for any function $f$ such that $f^{2}(\mu,\theta,\lambda) \le
V(\theta, \lambda)$ for all $(\mu,\theta,\lambda)$.  Of course, it is
easy to find functions involving $\mu$ or $\theta$ that do not satisfy
this requirement.  On the other hand, Theorem~\ref{thm:dm_clt} will be
useful for many functions of only $\lambda_{\theta}$ and
$\lambda_{e}$.

Recall that the drift function may not be unique.  Prior to the work
of \citet{jone:hobe:2004}, \citet{hobe:geye:1998} also analyzed this
Gibbs sampler and established \eqref{eq:geo_dc} using a different
drift function and more restrictive conditions on $a_{1}$ and $m'$.
However, this drift function can alleviate some of the difficulties
with using Theorem~\ref{thm:dm_clt} for functions involving $\mu$.

\begin{proposition}
\label{pr:hg_gibbs}
\citep{hobe:geye:1998}  Assume that $a_{1} \ge (3K-2)/(2K-2)$ and $m'
> (\sqrt{5} - 2) m''$. Then the Gibbs sampler satisfies
\eqref{eq:geo_dc} with the drift function 
\[
\begin{split}
  W(\mu, \theta, \lambda) = & 1 + \vartheta \left(
  \frac{1}{\lambda_{e}} +  \sum_{i=1}^{K} m_{i} (\bar{y}_{i} -
  \theta_{i})^{2} + (\mu - \bar{y})^{2} \right) + \\
& + \frac{1}{\lambda_{\theta}} + e^{c_1\lambda_{\theta}} +
  e^{c_1\lambda_{e}} + \frac{K\lambda_{\theta}}{s_{0} +
    K\lambda_{\theta}}(\overline{\theta} -\overline{y})^2\; .
\end{split}
\]
where $0 < \vartheta < 1$ is a constant defined on p. 427 in
\citet{hobe:geye:1998}.
\end{proposition}

Proposition~\ref{pr:jh_gibbs} shows that this Gibbs sampler is
geometrically ergodic as long as $a_1>3/2$ and $5m'>m''$.  However, it
does not satisfy detailed balance.  An appeal to
Corollary~\ref{cor:ge_pe} or ~\ref{cor:ge} shows that functions with a
little bit more than a second moment with respect to \eqref{eq:post}
will enjoy a CLT.

\subsection{Independence Sampler}
\label{sec:indep}
The independence sampler is an important special case of the MHG
algorithm.  Suppose the target distribution $\pi$ has support $\mathsf{X}
\subseteq \mathbb{R}^{k}$ and a density which, in a slight abuse of
notation, is also denoted $\pi$.  Let $p$ be a proposal density whose
support contains $\mathsf{X}$ and suppose the current state of the
chain is $X_n = x$.  Draw $y \sim p$ and set $X_{n+1}=y$ with
probability
\[
\alpha(x,y) = \frac{\pi(y) p(x)}{\pi(x) p(y)} \, \wedge \, 1 \; ;
\]
otherwise set $X_{n+1}=x$.  This Markov chain is Harris ergodic and it
is well-known \citep{meng:twee:1996} that it is uniformly ergodic if
there exists $\kappa > 0$ such that
\begin{equation}
\label{eq:indep}
\frac{\pi(x)}{p(x)} \le \kappa 
\end{equation}
for all $x$ since \eqref{eq:indep} implies a minorization
\eqref{eq:min_con} on $\mathsf{X}$ with $n_0 = 1$ and
$\epsilon=1/\kappa$.  Hence Corollary~\ref{cor:t94} implies that a CLT
will hold for $\bar{f}_{n}$ if $\text{E}_{\pi} f^{2} < \infty$.  On
the other hand, the independence sampler will not even be
geometrically ergodic if there is a set of positive $\pi$-measure
where \eqref{eq:indep} fails to hold.  Moreover, in this case
\citet{robe:1999} has given conditions which ensure a CLT
\textit{cannot} hold.

\subsection{An MHG Algorithm for Finite Point Processes}
\label{sec:fpp}
The material in this subsection is adapted from \citet{geye:1999} and
\citet{moll:1999}. Let ${\cal X}$ be a bounded region of
$\mathbb{R}^{d}$ and let $\lambda$ be Lebesgue measure.  Define ${\cal
  X}^{0} := \{ \varnothing\}$ and for $k \ge 1$, ${\cal X}^{k} :=
{\cal X} \times \cdots \times {\cal X}$ (there being $k$ terms in the
Cartesian product).  Think of $x \in {\cal X}^{k}$ as a pattern of $k$
points in ${\cal X}$, in particular, ${\cal X}^{0}$ denotes the
pattern with no points, and define $n(x)$ to be the cardinality of $x$
so that if $x \in {\cal X}^{k}$ then $n(x)=k$.  Let the state space
$\mathsf{X}$ be the union of all ${\cal X}^{k}$, that is, $\mathsf{X}=
\cup_{i=0}^{\infty} \mathsf{X}_{i}$ where $\mathsf{X}_{i}=\{x \, : \,
n(x)=i\}$. The target $\pi$ is an unnormalized density with respect
to the Poisson process with intensity measure $\lambda$ on ${\cal X}$.
\citet{geye:1999} proposes the following MHG algorithm for simulating
from $\pi$:

\vspace*{2mm}
\begin{enumerate}
\item With probability 1/2 attempt an up step
\begin{enumerate}
\item Draw $\xi \sim \lambda(\cdot)/\lambda({\cal X})$.  Set $x=x \cup
  \xi$ with probability
\[
1 \, \wedge \, \frac{\lambda({\cal X}) \, \pi(x \cup \xi)}{(n(x) + 1)
  \, \pi(x)} \; . 
\]
\end{enumerate}
\item Else attempt a down step
\begin{enumerate}
\item If $x = \varnothing$ skip the down step
\item Draw $\xi$ uniformly from the points of $x$.  Set $x=x
    \setminus \xi$ with probability
\[
1 \, \wedge \, \frac{n(x) \, \pi(x \setminus \xi)}{\lambda({\cal X})
  \, \pi(x)} \; . 
\]
\end{enumerate}
\vspace*{-3.5mm}
\end{enumerate}

This MHG algorithm is Harris recurrent and geometrically ergodic.
\begin{proposition} \citep{geye:1999}
\label{pr:geye_pp}
Suppose there exists a real number $M$ such that 
\[
\pi(x \cup \xi) \le M \pi(x)
\]
for all $x \in \mathsf{X}$ and all $\xi \in {\cal X}$.  Then the MHG
algorithm started at $x^* \in \{ x \, : \, \pi(x) >0\}$ is Harris
ergodic and satisfies \eqref{eq:geo_dc} with the drift function $V(x)
= A^{n(x)}$ where $A > M \lambda({\cal X})\, \vee\, 1$.
\end{proposition}

Of course, Theorem~\ref{thm:dm_clt} implies a CLT for $\bar{f}_{n}$
for any function $f$ such that $f^{2}(x) \le A^{n(x)}$ for all $x$.  On
the other hand, this algorithm was constructed so as to satisfy
\eqref{eq:dbc} (see \citet{geye:1999} for a detailed argument) and
hence the Markov chain is asymptotically uncorrelated so that a CLT
holds when $\text{E}_{\pi} f^{2}(x) < \infty$.

\subsection{Random Walk MHG Algorithms}
\label{sec:rwmh}
Let $\pi$ be a target density on $\mathbb{R}^{k}$ and let the proposal
density have the form $q(y|x)=q(|y-x|)$.  Now suppose that the current
state of the chain is $X_{n}=x$.  Draw $y \sim q$ and set $X_{n+1}=y$
with probability
\[
\alpha(x,y) = \frac{\pi(y)}{\pi(x)} \, \wedge \, 1 \; ;
\]
otherwise set $X_{n+1}=x$.  Note that this algorithm satisfies
\eqref{eq:dbc} by construction.

Random walk-type MHG algorithms are some of the most useful and
popular MCMC algorithms and consequently their theoretical properties
have been thoroughly studied. \citet{meng:twee:1996} show that random
walk samplers (on $\mathbb{R}^{k}$) \textit{cannot} be uniformly
ergodic (or uniformly mixing) but they do establish that a random walk
MHG algorithm can be geometrically ergodic by verifying
\eqref{eq:geo_dc} when $k=1$ and $\pi$ has tails that decrease
exponentially. \citet{robe:twee:1996} extended their work by
establishing \eqref{eq:geo_dc} in the case where $k \ge 1$.  However,
\citet{jarn:hans:2000} verified \eqref{eq:geo_dc} with a different
drift function than that used by \citet{robe:twee:1996} and obtained
more general conditions ensuring geometric ergodicity.  On the other
hand, if a random walk MHG algorithm is not geometrically ergodic it
may still be polynomially ergodic of all orders; see
\citet{fort:moul:2000}.

\begin{proposition}
\label{pr:rwmhg}
\citep{jarn:hans:2000} Suppose $\pi$ is a positive density on
$\mathbb{R}^{k}$ having continuous first derivatives such that
\[
\lim_{|x| \rightarrow \infty} \frac{x}{|x|} \cdot \nabla \log \pi(x) =
- \infty \; .
\] 
Let $A(x) := \{ y \in \mathbb{R}^{k} \, : \, \pi(y) \ge \pi(x) \}$ be
the region of certain acceptance and assume that there exist
$\delta>0$ and $\epsilon >0$ such that, for every $x$, $|x-y| \le
\delta$ implies $q(y|x) \ge \epsilon$.  Then if
\[
\liminf_{|x| \rightarrow \infty} \int_{A(x)} q(y|x)\,dy > 0
\]
the random walk MHG algorithm satisfies \eqref{eq:geo_dc} with the
drift function $V(x) = c\pi(x)^{-1/2}$ for some $c > 0$.
\end{proposition}

Hence, under the conditions of Proposition~\ref{pr:rwmhg},
Theorem~\ref{thm:dm_clt} guarantees a CLT if $f(x)$ satisfies
$f^{2}(x) \le c\pi(x)^{-1/2}$ for all $x \in \mathbb{R}^{k}$.
Alternatively, we conclude that the random walk MHG is geometrically
ergodic, satisfies \eqref{eq:dbc} (and hence is asymptotically
uncorrelated) and an appeal to Corollary~\ref{cor:rr97} establishes
the existence of a CLT if $\text{E}_{\pi} f^{2}(x) < \infty$.

\section{Final Remarks}
\label{sec:remarks}
The focus has been on some of the connections between recent work on
general state space Markov chains and results from mixing processes
and the implications for Markov chain CLTs.  However, this article
only scratches the surface of the mixing process literature that is
potentially useful in MCMC.  For example, the existence of a
functional CLT or strong invariance principle is required in order to
estimate $\sigma^{2}_{f}$ from \eqref{eq:clt}
\citep{dame:1994,glyn:whit:1992,jone:hara:caff:2004}.  There has been
much work on these for mixing processes; \citet{phil:stou:1975} is a
good starting place for strong invariance principles while
\citet{bill:1968} gives an introduction to the functional CLT.

\begin{appendix}
\section{Calculations for Example~\ref{ex:dmc}} 
\label{app:dmc}
Define $V(z) =a^{|z|}$ for some $a > 1$.  Then $V(z) \ge 1$ for all $z
\in \mathbb{Z}$ and
\[
PV(x) = \sum_{y \in \mathbb{Z}} a^{|y|} P(x,y) \; .
\] 
Recall that $\Delta V(x) = PV(x) - V(x)$.  The first goal is to show
that if $x \neq 0$ then $ \Delta V(x)/V(x)<0$ since then there must be
a $\beta > 0$ such that $ \Delta V(x)< -\beta V(x)$.  Suppose $X_{n}=x
\ge 1$ then
\[
PV(x) =  \theta a^{x+1} + 1-\theta \; \Rightarrow \; \frac{\Delta
  V(x)}{V(x)} = a \theta - 1 + \frac{1 - \theta}{a^{x}} < 0 
\]
as long as
\begin{equation}
\label{eq:aa}
(a \theta - 1)V(x) + 1 - \theta < 0 \; .
\end{equation}
Now \eqref{eq:aa} can hold only if $ a\theta - 1 < 0$ and since $V(x)
\ge a$ for all $x \neq 0$ \eqref{eq:aa} will hold when $ a \theta - 1
< 0$ and $(a \theta - 1)a + 1 - \theta < 0$.  A similar argument shows
that this is also the case when $X_{n}=-x \le -1$.  Now suppose
$X_{n}=0$. Then $PV(0)=a$ and $\Delta V(0) = - V(0) + a$.  Putting
this together yields \eqref{eq:geo_dc} with $C=\{0\}$.

\end{appendix}

\bigskip
\bigskip
\noindent {\bf Acknowledgments}
\bigskip

\noindent The author is grateful to Charlie Geyer, Murali Haran, Jeff
Rosenthal and Qi-Man Shao for helpful conversations about this paper.
Also, the author thanks an anonymous referee for many very helpful
comments.

\end{document}